\documentclass{svproc}
\usepackage{amsmath, amssymb}       
\usepackage{verbatimbox}
\usepackage{verbatim}
\usepackage{tcolorbox}

\usepackage{url}

\begin{document}
\mainmatter              
\title{Numerical Continuation and Bifurcation Analysis in a Harvested Predator-Prey Model with Time Delay using DDE-Biftool}
\titlerunning{Numerical Continuation and Bifurcation Analysis using DDE-Biftool}  
%
\author{Juancho A. Collera}
\authorrunning{Juancho A. Collera} 
%
\tocauthor{Juancho A. Collera}
\institute{University of the Philippines Baguio,\\ Gov. Pack Road, Baguio City 2600 Philippines\\
\email{jacollera@up.edu.ph}}

\maketitle              

\begin{abstract}
Time delay has been incorporated in models to reflect certain physical or biological meaning. 
The theory of delay differential equations (DDEs), which has seen extensive growth in the last seventy years or so, can be used to examine the effects of time delay in the dynamical behavior of systems being considered. 
Numerical tools to study DDEs have played a significant role not only in illustrating theoretical results but also in discovering interesting dynamics of the model. DDE-Biftool, which is a Matlab package for numerical continuation and numerical bifurcation analysis of DDEs, is one of the most utilized and popular numerical tools for DDEs. 
In this paper, we present a guide to using the latest version of DDE-Biftool targeted to researchers who are new to the study of time delay systems. 
A short discussion of an example application, which is a harvested predator-prey model with a single discrete time delay, will be presented first. 
We then implement this example model in DDE-Biftool, pointing out features where beginners need to be cautious. 
We end with a comparison of our theoretical and numerical results.
\end{abstract}
\section{Introduction}
Time delay has been incorporated in models to reflect certain physical or biological meaning. 
Examples include optical feedback in laser systems \cite{BuonoCollera2015,Collera2015,Collera2016}, maturation age in stage structured population models \cite{ColleraMagpantay2018}, and delayed information in queueing models \cite{Penderetal2017} just to name a few. The theory of delay differential equations (DDEs) \cite{HaleLunel1993,Smith2011}, which has seen extensive growth in the last seventy years or so, can be used to examine the effects of time delay in the dynamical behavior of systems being considered. Numerical tools to study DDEs have played a significant role not only in illustrating theoretical results but also in discovering interesting dynamics of the model. DDE-Biftool \cite{Engelborghsetal2001,Sieberetal2016}, which is a Matlab package for numerical continuation and numerical bifurcation analysis of DDEs, is one of the most utilized and popular numerical tools for DDEs. \\

In this paper, we present a guide to using the latest version of DDE-Biftool targeted to researchers who are new to the study of time delay systems. A short theoretical discussion of the example, which is a harvested predator-prey model with a single discrete time delay, will be given in the next section.  In section 3, we implement this example model in DDE-Biftool and compare the theoretical and numerical results. We conclude the paper with a summary and thoughts on using DDE-Biftool in studying time delay systems.

\section{Harvested Predator-Prey Model with Time Delay}
We first discuss the model that we are going to use as an example for numerical continuation and numerical bifurcation analysis. 
The theoretical results presented here will be compared to the numerical results obtained in the succeeding section.

\subsection{The Model}

We consider the following model, 
studied in \cite{ToahaHassan2008}, with a single discrete time delay parameter $\tau>0$ 

\begin{equation}\label{eq:model}
\left\{
\begin{tabular}{ccl}
$\displaystyle \frac{dx(t)}{dt}$ 
& $\ \ =\ \ $ & $rx(t)-ax(t) x(t-\tau) - bx(t)y(t) - h$,\\
&&\\
$\displaystyle \frac{dy(t)}{dt}$ 
& $=$ & $cx(t)y(t) - dy(t) - k$.\\
\end{tabular}
\right.
\end{equation}
Here, $x(t)$ and $y(t)$ are the state variables denoting, respectively, the densities of the prey and predator populations at time $t$. 
The parameter $a$ is the ratio of the intrinsic growth rate $r$ of the prey and the carrying capacity $K$ for the prey population in the absence of the predation. 
The rate of consumption of prey by the predator is given by the parameter $b$ while $c$ measures the conversion of prey consumed into the predator reproduction rate. 
The death rate of the predator is represented by $d$. 
Both species are assumed to have economic value and are harvested. The parameters $h$ and $k$ denote the harvesting rates of the prey and predator populations, respectively. 
All parameters in system (\ref{eq:model}) are positive real numbers.\\ 

The equilibrium solutions of system (\ref{eq:model}) are solutions $(x(t),y(t))$ of system (\ref{eq:model}) satisfying 
$dx(t)/dt=0$ and $dy(t)/dt=0$, 
and hence are obtained by solving for constant values $x$ and $y$ in the following system of nonlinear equations

\begin{equation} \label{eq:equations}
\left\{
\begin{tabular}{rcl}
$rx-ax^2-bxy-h$ & $\ \ =\ \ $ & $0$,\\
&&\\
$cxy-dy-k$ & $=$ & $0$.\\
\end{tabular}
\right.
\end{equation}
Since $h>0$ and $k>0$, the value of $x$ nor the value of $y$ can not be zero. 
From the equations in system (\ref{eq:equations}), we get $y=(rx-ax^2-h)/bx$ and $y=k/(cx-d)$. Hence, the positive equilibria of system (\ref{eq:model}) exist provided $(r^2-4ah)>0$ and $(cx-d)>0$.

\begin{table}[h]
\begin{center}
\begin{tabular}{|c|c|c|c|c|c|c|} \hline
$r$ & $a$ & $b$ & $c$ & $d$ & $h$ & $k$\\ \hline
\ \ \ 3.50\ \ \  & 
\ \ \ 0.04\ \ \  &  
\ \ \ 1.00\ \ \  & 
\ \ \ 0.05\ \ \  & 
\ \ \ 0.30\ \ \  & 
\ \ \ 0.02\ \ \  & 
\ \ \ 0.01\ \ \ \\ \hline
\end{tabular}
\caption{Parameter values.}
\label{tab:param}
\end{center}
\end{table}

\begin{example}
Solving the nonlinear equations in system (\ref{eq:equations}) with parameter values given in Table \ref{tab:param}, we obtain the following equilibrium solutions of system (\ref{eq:model})
\begin{eqnarray} 
E_1 &=& (\ \ 6.061458,\ \  3.254242), \label{eq:e1}\\
E_2 &=& (87.432881,\ \  0.002456). \label{eq:e2}
\end{eqnarray}
\end{example}

\subsection{Local Stability of the Equilibrium Solutions}
If we let $\overline{x}(t)=x(t-\tau)$ and $\overline{y}(t)=y(t-\tau)$, 
then the right-hand side of system (\ref{eq:model}) can be expressed as
\begin{equation}\label{eq:model2}
\left[
\begin{tabular}{c}
$f(x,y,\overline{x},\overline{y})$\\
$g(x,y,\overline{x},\overline{y})$\\
\end{tabular}
\right]
=
\left[
\begin{tabular}{c}
$rx-ax\overline{x} - bxy - h$\\
$cxy - dy - k$\\
\end{tabular}
\right].
\end{equation}
The linearized system corresponding to system (\ref{eq:model}) about an equilibrium solution $E^*=(x^*,y^*)$ is given by 
\begin{equation}\label{eq:linsys}
\dfrac{d\mathbf{X}(t)}{dt} = \mathbf{A} \mathbf{X}(t) + \mathbf{B} \mathbf{X}(t-\tau)
\end{equation}
where 
$
\mathbf{X}(t) 
=
\left[
\begin{tabular}{c}
$x(t)$\\
$y(t)$\\
\end{tabular}
\right]
$ 
and the matrices $\mathbf{A}$ and $\mathbf{B}$ are as follows 
\begin{eqnarray}
\mathbf{A} 
&=&
\left[
\begin{tabular}{ccc}
$\partial f/\partial x$ & $\quad$ & $\partial f/\partial y$\\
$\partial g/\partial x$ & $\quad$ & $\partial g/\partial y$\\
\end{tabular}
\right]
_{(x,y,\overline{x},\overline{y})=(x^*,y^*,x^*,y^*)}\label{eq:A}
,\\
\mathbf{B} 
&=&
\left[
\begin{tabular}{ccc}
$\partial f/\partial \overline{x}$ & $\quad$ & $\partial f/\partial \overline{y}$\\
$\partial g/\partial \overline{x}$ & $\quad$ & $\partial g/\partial \overline{y}$\\
\end{tabular}
\right]
_{(x,y,\overline{x},\overline{y})=(x^*,y^*,x^*,y^*)}\label{eq:B}
.
\end{eqnarray}
The characteristic equation corresponding to the linear system (\ref{eq:linsys}) is 
\begin{equation}\label{eq:ce0}
\det(\lambda \mathbf{I} - \mathbf{A} - \mathbf{B} e^{-\lambda\tau})=0
\end{equation}
where $\mathbf{I}$ is the $2\times 2$ identity matrix. 
Equation (\ref{eq:ce0}) is obtained by using the ansatz 
$
\mathbf{X}(t) 
=
e^{\lambda t}\left[
\begin{tabular}{c}
$x_0$\\
$y_0$\\
\end{tabular}
\right]
$ 
to the linear system (\ref{eq:linsys}). 
If all roots of the characteristic equation (\ref{eq:ce0}) lie in the open left-half plane, i.e. $\mbox{Re}\ \lambda <0$ for all roots $\lambda$ of equation (\ref{eq:ce0}), then the equilibrium $E^*$ is locally asymptotically stable.  \\

Using equations (\ref{eq:A}) and (\ref{eq:B}) and the functions $f$ and $g$ in (\ref{eq:model2}), we obtain 
$$
\mathbf{A} 
=
\left[
\begin{tabular}{ccc}
$r-ax^*-by^*$ & $\quad$ & $-bx^*$\\
$cy^*$ & $\quad$ & $cx^*-d$\\
\end{tabular}
\right]
\quad\mbox{and}\quad
\mathbf{B} 
=
\left[
\begin{tabular}{ccc}
$-ax^*$ & $\quad$ & $0$\\
$0$ & $\quad$ & $0$\\
\end{tabular}
\right].
$$
Thus, we can write the characteristic equation (\ref{eq:ce0}) as 
\begin{equation}\label{eq:ce}
\left( \lambda^2+a_1\lambda+a_2 \right) 
+ \left( a_3\lambda+a_4 \right) e^{-\lambda\tau} \ \ =\ \ 0
\end{equation}
where 
$a_1 = -(r-ax^*-by^*) -(cx^*-d)$, 
$a_2 = (r-ax^*-by^*)(cx^*-d) + bcx^*y^*$, 
$a_3 = ax^*$, and 
$a_4 = -ax^*(cx^*-d)$. \\

When the time delay $\tau=0$, equation (\ref{eq:ce}) reduces to the quadratic equation  
\begin{equation}{\label{eq:poly}}
\lambda^2+(a_1+a_3)\lambda+(a_2+a_4)\ \ =\ \ 0.
\end{equation}
Both roots of equation (\ref{eq:poly}) have negative real part if and only if 
\begin{equation}{\label{eq:polycond}}
(a_1+a_3)>0
\qquad\mbox{and}\qquad
(a_2+a_4)>0.
\end{equation}
Hence, for the case $\tau=0$, the equilibrium $E^*$ is locally asymptotically stable whenever conditions in (\ref{eq:polycond}) are satisfied.\\
 
We wanted to know if $E^*$, under the conditions in (\ref{eq:polycond}) will become unstable as we vary the time delay parameter. Suppose that the conditions in (\ref{eq:polycond}) are satisfied and consider now the case where $\tau>0$. Initially, $E^*$ is locally asymptotically stable, i.e. all roots of the characteristic equation (\ref{eq:ce}) with $\tau=0$ lie in the open left-half plane. If one or more roots of equation (\ref{eq:ce}) cross the imaginary axis and move towards the open right-half plane as $\tau$ is increased, then $E^*$ will switch stability and becomes unstable. We have two possibilities: either a real root of equation (\ref{eq:ce}) will cross the imaginary axis, i.e. $\lambda=0$ is a root of equation (\ref{eq:ce}) at some critical delay value, or a pair of complex conjugate roots of equation (\ref{eq:ce}) cross the imaginary axis, i.e. $\lambda=\pm i\omega$ is a root of equation (\ref{eq:ce}) at some critical delay value where $\omega$ is a nonzero real number.\\

If $\lambda=0$ is a root of equation (\ref{eq:ce}), then $(a_2+a_4)=0$. 
However, since $(a_2+a_4)>0$ from conditions in (\ref{eq:polycond}), then $\lambda=0$ is not a root of the characteristic equation (\ref{eq:ce}). 
Suppose now that equation (\ref{eq:ce}) has a pair of purely imaginary roots $\lambda=\pm i\omega$. 
Since the right-hand side of equation (\ref{eq:ce}) is an entire function, complex roots of equation (\ref{eq:ce}) come in conjugate pairs. 
Thus, without loss of generality, we may assume that $\omega>0$. 
Since $\lambda=i\omega$ with $\omega>0$ satisfies equation (\ref{eq:ce}), we have 
\begin{equation}\label{eq:sub}
(-\omega^2+ia_1\omega+a_2)+(ia_3\omega+a_4)e^{-i\omega\tau} = 0.
\end{equation} 
This gives the following equations 
\begin{eqnarray}
a_4\cos (\omega\tau) + a_3\omega\sin (\omega\tau) & \ \ =\ \ & \omega^2-a_2, \label{eq:real}\\
a_4\sin (\omega\tau) - a_3\omega\cos (\omega\tau) & = & a_1\omega, \label{eq:imag}
\end{eqnarray}
after using the Euler's formula in equation (\ref{eq:sub}) and then matching the real and imaginary parts on both sides of equation (\ref{eq:sub}). 
We can eliminate $\tau$ by squaring each side of equations (\ref{eq:real}) and (\ref{eq:imag}) and then adding corresponding sides. 
We obtain 
$
(\omega^2-a_2)^2 + a_1^2\omega^2 = a_4^2 + a_3^2\omega^2, 
$ 
which we can write as 
\begin{equation} \label{eq:deg4}
\omega^4 + \alpha\omega^2 + \beta\ \ =\ \ 0
\end{equation}
where
$\alpha =  a_1^2-2a_2-a_3^2$ and $\beta =  a_2^2-a_4^2$. 
If we let $u=\omega^2$, then equation (\ref{eq:deg4}) becomes the following quadratic equation in $u$
\begin{equation} \label{eq:deg2}
h(u)\ \ :=\ \ u^2+\alpha u+\beta\ \ =\ \ 0.
\end{equation}

If equation (\ref{eq:deg2}) does not have a positive root, 
then equation (\ref{eq:ce}) cannot have purely imaginary roots. 
That is, the roots of the charcateristic equation (\ref{eq:ce}) that are in the open left-half plane when $\tau=0$ remain in the open left-half plane as the time delay parameter $\tau$ is increased. 
In other words, if equation (\ref{eq:deg2}) does not have a positive root, then the equilibrum $E^*$ remains localy asymptotically stable for all $\tau>0$. 
Note that if the coefficients in equation (\ref{eq:deg2}) satify the following conditions
\begin{equation} \label{eq:absstab}
\alpha > 0 
\qquad \mbox{and} \qquad
\beta > 0,
\end{equation}
then both roots of equation (\ref{eq:deg2}) have negative real parts. 
That is, under the conditions in (\ref{eq:absstab}), equation (\ref{eq:deg2}) does not have positive roots.  
Therefore, the equilibrium $E^*$ is locally asymptotically stable for all $\tau\ge0$ whenever conditions in (\ref{eq:polycond}) and (\ref{eq:absstab}) are satisfied (cf. Theorem 4 of \cite{ToahaHassan2008}).

\begin{lemma}{\label{lem:posroots}}
The number of positive roots of equation (\ref{eq:deg2}) is determined as follows 
\begin{center}
\begin{tabular}{|c|c|}\hline
 Conditions & \ \ \ \ \ Number of Positive Roots\ \ \ \ \ \\ \hline
\ \ \ \ \ $\alpha\ge0$ \ \ and\ \ \  $\beta\ge0$\ \ \ \ \  & $0$\\ \hline
$\alpha\ge0$ \ \  and\ \  \ $\beta<0$ & $1$\\ \hline
$\alpha<0$ \ \ and\ \ \  $\beta\le0$ & $1$\\ \hline 
$\alpha<0$ \ \ and\ \ \  $\beta>0$\ \ \ with\ \ \ $h(\overline{u})>0$   & $0$\\ \hline 
$\alpha<0$ \ \ and\ \ \  $\beta>0$\ \ \ with\ \ \ $h(\overline{u})=0$   & $1$\\ \hline 
$\alpha<0$ \ \ and\ \ \  $\beta>0$\ \ \ with\ \ \ $h(\overline{u})<0$   & $2$\\ \hline 
\end{tabular}
\end{center}
where $(\overline{u},h(\overline{u}))$ is the vertex of the parabola given by the graph of the function $h(u)$ in equation (\ref{eq:deg2}).
\end{lemma}

\begin{example}\label{ex:two}
Using the parameter values in Table \ref{tab:param}, 
we obtain $\alpha = -2.031311$ and $\beta = 0.972753$ approximately. 
The graph of the quadratic function $h(u)$ given in equation (\ref{eq:deg2}) is a parabola with vertex at $(\overline{u},h(\overline{u}))=(1.015655, -0.058803)$. 
Since $\alpha<0$, $\beta>0$, and $h(\overline{u})<0$, by Lemma \ref{lem:posroots}, equation (\ref{eq:deg2}) has two positive roots. 
Let us denote these positive roots by $u_-$ and $u_+$ with $u_- < u_+$. 
Solving equation (\ref{eq:deg2}), we obtain the roots $u_- = 0.773162$ and $u_+ = 1.258149$, with corresponding 
\begin{equation} \label{eq:omegas}
\omega_- = 0.879297
\qquad \mbox{and} \qquad
\omega_+ = 1.121672,
\end{equation}
which are roots of equation (\ref{eq:deg4}). 
Consequently, the characteristic equation (\ref{eq:ce}) has purely imaginary roots $\pm i\omega_-$ and $\pm i\omega_+$.
\end{example}

\begin{remark}
The conditions in (3.13) of \cite{ToahaHassan2008} are in fact the three inequalities in the last row of the table in Lemma \ref{lem:posroots}, i.e.

\begin{equation} \label{eq:cond2posroots}
\alpha < 0, \qquad 
\beta>0,  
\qquad \mbox{and} \qquad
h(\overline{u})<0,
\end{equation}
where equation (\ref{eq:deg2}) has exactly two positive roots.
\end{remark}

\subsection{Critical Delay Values}
Let us now determine the critical time delay values where the purely imaginary roots $\pm i\omega_-$ and $\pm i\omega_+$ of equation (\ref{eq:ce}), obtained in Example \ref{ex:two}, will occur. 
Under the conditions in (\ref{eq:cond2posroots}), equation (\ref{eq:ce}) has purely imaginary roots $\pm i\omega_{-}$ and $\pm i\omega_{+}$. 
Thus, the values $\pm\omega_-$ and $\pm\omega_+$ satisfy equations (\ref{eq:real}) and (\ref{eq:imag}). 
We can compute for $\sin(\omega\tau)$ and $\cos (\omega\tau)$ from equations (\ref{eq:real}) and (\ref{eq:imag}) to obtain 
$$ 
\tan(\omega\tau)
=
\frac{\omega(a_3\omega^2+a_1a_4-a_2a_3)}{(a_4-a_1a_3)\omega^2-a_2a_4}.
$$ 
The purely imaginary roots $\pm i\omega_-$ (resp. $\pm i\omega_+$) of equation (\ref{eq:ce}) occurs when the time delay $\tau=\tau_k^-$ (resp. $\tau=\tau_k^+$) for $k = 0,1,2,3,\dots$ with 
\begin{equation} \label{eq:tauk}
\tau_k^\pm
=
\frac{1}{\omega_\pm}
\left[
\tan^{-1}\left(\frac{\omega_\pm(a_3\omega_\pm^2+a_1a_4-a_2a_3)}{(a_4-a_1a_3)\omega_\pm^2-a_2a_4}\right) + 2\pi k
\right].
\end{equation}
\begin{example}\label{ex:three}
Using equation (\ref{eq:tauk}) with the parameter values in Table \ref{tab:param} and corresponding $\omega_\pm$ given in (\ref{eq:omegas}), we get the following critical time delay values. 
\begin{table}
\centering
\begin{tabular}{|lcr||lcr|} \hline
\ \ \ \ $\tau_0^-$ & $=$ & $-1.752556$\ \ \ \ & \ \ \ \ $\tau_0^+ $ & $=$ &   $1.3794139$\ \ \ \ \\ \hline
\ \ \ \ $\tau_1^- $ & $=$ &   $5.393140$\ \ \ \  & \ \ \ \ $\tau_1^+ $ & $=$ &   $6.9810371$\ \ \ \ \\ \hline 
\ \ \ \ $\tau_2^- $ & $=$ & $12.538836$\ \ \ \  & \ \ \ \ $\tau_2^+ $ & $=$ & $12.582660$\ \ \ \ \\ \hline
\ \ \ \ $\tau_3^- $ & $=$ & $19.684531$\ \ \ \  & \ \ \ \ $\tau_3^+ $ & $=$ & $18.184284$\ \ \ \ \\ \hline
\ \ \ \ $\tau_4^- $ & $=$ & $26.830227$\ \ \ \  & \ \ \ \ $\tau_4^+ $ & $=$ & $23.785907$\ \ \ \ \\ \hline
\end{tabular}
\caption{Values of $\tau_k^-$ and $\tau_k^+$ for $k=0,1,2,3,4$.}
\label{tab:taus}
\end{table}
\end{example}

\subsection{Transversality Conditions}
We saw earlier that at $\tau=\tau_k^-$ (resp. $\tau=\tau_k^+$) for $k = 0,1,2,3,\dots$, the characteristic equation (\ref{eq:ce}) has purely imaginary roots $\pm i\omega_-$ (resp. $\pm i\omega_+$). 
We wanted to know if these roots along the imaginary axis will move towards the open right-half plane or towards the open left-half plane. 
We address this by determining if the rate of change of the $\mbox{Re}\ \lambda$ with respect to $\tau$ at the critical time delay values $\tau_\pm$ is positive or negative, where $\lambda=\lambda(\tau)$ is a root of the characteristic equation (\ref{eq:ce}). \\

Recall, from equation (\ref{eq:deg2}), that $h(u)=u^2+\alpha u+\beta$. 
Hence,  $h'(u) = 2u + \alpha$. 
As shown in equation (3.15) of \cite{ToahaHassan2008}, we have 
$$
\left.\mbox{sign}\left\{\frac{d(\mbox{Re}\lambda)}{d\tau}\right\}\right|_{
\tau\ =\ \tau_k^\pm }
\ =\ \ 
\mbox{sign}\left\{h'(\omega_\pm^2)\right\}
\ =\ \ 
\mbox{sign}\left\{h'(u_\pm)\right\}
$$
in our notation. Since the graph of $h(u)$ is decreasing (resp. increasing) at $u=u_-$ (resp. at $u=u_+$), we know that $h'(u_-)<0$ (resp. $h'(u_+)>0$). Therefore, 
$$
\left.\mbox{sign}\left\{\frac{d(\mbox{Re}\lambda)}{d\tau}\right\}\right|_{
\tau\ =\ \tau_k^\pm }
\ =\ \ 
\pm 1.
$$
This means that the root $\lambda(\tau)$ of the characteristic equation (\ref{eq:ce}) that lies on the imaginary axis when $\tau=\tau_k^-$ (resp. when $\tau=\tau_k^+$) moves towards the open left-half plane (resp. towards the open right-half plane).

\section{Numerical Continuation and Bifurcation Analysis}
Numerical tools to study DDEs have played a significant role not only in illustrating theoretical results but also in discovering interesting dynamics of the model. 
\emph{DDE-Biftool} \cite{Engelborghsetal2001,Sieberetal2016}, which is a Matlab package for numerical continuation and numerical bifurcation analysis of DDEs, is one of the most utilized and popular numerical tools for DDEs. It was originally developed by K. Engelborghs \cite{Engelborghsetal2001} as part of his PhD work at the KU Leuven under supervision of D. Roose. 
DDE-Biftool provides a set of capabilities that is similar to what a range of alternative tools do for ordinary differential equations (ODEs) and maps, such as \emph{Matcont} \cite{Dhoogeetal2003}, \emph{COCO} \cite{DankowiczSchilder2013} and \emph{AUTO} \cite{DoedelOldeman2012}. 
Another tool performing a similar set of tasks for DDEs, particularly for time-dependent DDEs with time-dependent delays, is \emph{Knut} \cite{Szalai2013}.\\

Aside from continuation of steady-state and periodic-orbit solutions which are typically done by varying a single parameter, DDE-Biftool can also continue bifurcations in two parameters. This includes steady-state folds, Hopf bifurcations, folds of periodic orbits, period doublings, and torus bifurcations. It can also perform normal form analysis for equilibria. 
DDE-Biftool is GNU Octave compatible and has a BSD licence such that it can be run completely as free software. 
The most recent version, DDE-Biftool v3.1.1, is maintained by J. Sieber and can be downloaded from https://sourceforge.net/projects/ddebiftool. The manual for this newest version of DDE-Biftool is provided at \cite{Sieberetal2016}.\\

In this section, we illustrate the use of DDE-Biftool to perform numerical continuation and bifurcation analysis of system (\ref{eq:model}) varying the delay parameter $\tau$. The boxed commands are the required commands and can be saved in a single m-file for convenience. \\

\begin{tcolorbox}
\begin{scriptsize}
\begin{verbatim}
addpath('../ddebiftool/','../ddebiftool_utilities/');

pp_sys = @(x,p)[...
p(1)*x(1,1) - p(2)*x(1,1).*x(1,2) - p(3)*x(1,1).*x(2,1) - p(6);
p(4)*x(1,1).*x(2,1) - p(5)*x(2,1)-p(7)];

funcs = set_funcs('sys_rhs',pp_sys,'sys_tau',@()[8]);
\end{verbatim}
\end{scriptsize}
\end{tcolorbox}

We start with \texttt{addpath} which identifies the location of the folders \texttt{ddebiftool} and \texttt{ddebiftool$\_$utilities} containing the functions that we need for continuation and bifurcation analysis. This should be adjusted depending on where the user intends to do their computations and where the folders  \texttt{ddebiftool} and \texttt{ddebiftool$\_$utilities} were copied.\\

Next, we encode the right-hand side of system (\ref{eq:model}) in the function named \texttt{pp$\_$sys}. 
The state variables $x_1(t)$, $x_1(t-\tau)$, $x_2(t)$, and $x_2(t-\tau)$ are respectively denoted by \texttt{x(1,1)}, \texttt{x(1,2)}, \texttt{x(2,1)}, and \texttt{x(2,2)}. Here, the first index refers to the component while the second index refers to the delay number. For simplicity, we represent the parameters $r$, $a$, $b$, $c$, $d$, $h$, and $k$ by \texttt{p(1)}, \texttt{p(2)}, \texttt{p(3)}, \texttt{p(4)}, \texttt{p(5)}, \texttt{p(6)}, and \texttt{p(7)}, respectively. \\

We then set-up the function structure \texttt{funcs} identifying the previously defined system right-hand side \texttt{pp$\_$sys} as \texttt{'sys$\_$rhs'}. 
The time delay $\tau$ will be the eighth parameter in our parameter list. 
Thus, assigning \texttt{[8]} in \texttt{'sys$\_$tau'}.\\

\begin{tcolorbox}
\begin{scriptsize}
\begin{verbatim}
parbd = {'min_bound',[8,0],'max_bound',[8,15],'max_step',[8,0.05]};

[br,success] = SetupStst(funcs,...
    'parameter',[3.50 0.04 1.00 0.05 0.30 0.02 0.01 0.00],...
    'x',[6.00; 3.00],'contpar',8,'step',0.02, parbd{:})    
\end{verbatim}
\end{scriptsize}
\end{tcolorbox}

In \texttt{parbd}, we set the minimum bound, maximum bound, and maximum stepsize for the time delay $\tau$, which is our main continuation and bifurcation parameter. 
Here, we choose $\tau$ to be from $0$ to $15$ because $\tau$ is non-negative and since we want to see the dynamics as we vary the time delay $\tau$ up until $\tau_2^+ = 12.582660$ (see Table \ref{tab:taus}). \\

Next, we set-up the branch of equilibria which we denote here by \texttt{br}. 
For \texttt{'parameter'}, we use the parameter values from Table \ref{tab:param}. 
Note that the intial value set for $\tau$, which is eighth in the parameter list, is \texttt{0.00}. 
For \texttt{'x'}, we use the initial guess \texttt{[6.00; 3.00]} targeting the equilibrium $E_1$ with values given in equation (\ref{eq:e1}). 
Our continuation parameter is $\tau$, so \texttt{'contpar'} is \texttt{8}. 
After running the commands above, we get the following results. 
\begin{scriptsize}
\begin{verbatim}
br = method: [1x1 struct]
     parameter: [1x1 struct]
     point: [1x2 struct]

success = 1
\end{verbatim}
\end{scriptsize}
This means that our attempt to set-up a branch of equilibria is successful. 
The equilibrium branch \texttt{br} now contains two points. 
The first branch point in \texttt{br} has $\tau=0$ while the second branch point in \texttt{br} has $\tau=0.02$ since \texttt{'step'} is assigned a value \texttt{0.02}. 
The corrected value for $E_1$ can be obtained by typing the following in the Command Window.
\begin{scriptsize}
\begin{verbatim}
>> format long; br.point(1).x

ans = 6.061458241811056
      3.254242134274988
\end{verbatim}
\end{scriptsize}
The value of the time delay parameter $\tau$ for the first and second points in the equilibrium branch \texttt{br} are obtained by typing the following in the Command Window.
\begin{scriptsize}
\begin{verbatim}
>> br.point(1).parameter(8)

ans = 0
\end{verbatim}
\end{scriptsize}
\begin{scriptsize}
\begin{verbatim}
>> br.point(2).parameter(8)

ans = 0.020000000000000
\end{verbatim}
\end{scriptsize}

\begin{figure}
\includegraphics[width=0.5\textwidth]{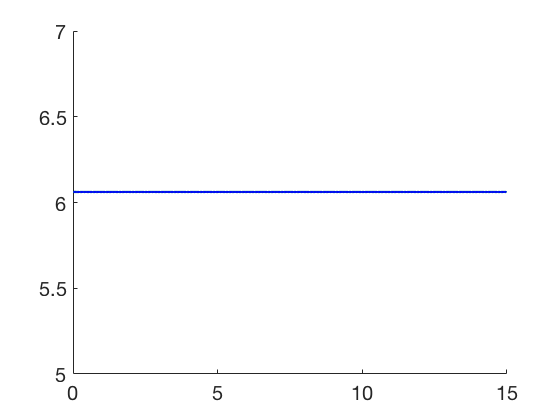}
\includegraphics[width=0.5\textwidth]{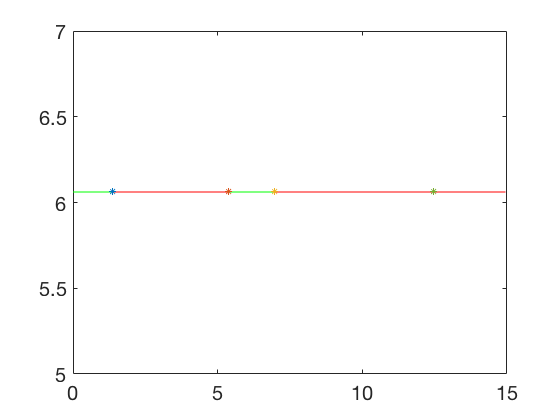}
\caption{(Left) The equilibrium branch \texttt{br} obtained by using the branch continuation function \texttt{br$\_$contn}. (Right) The same branch \texttt{br} with stability information obtained by using the function \texttt{br$\_$stabl}.}
\label{fig:contnsplot}
\end{figure}

\begin{tcolorbox}
\begin{scriptsize}
\begin{verbatim}
figure(1); clf;
br.method.continuation.plot = 1; 
[br,s,f,r] = br_contn(funcs,br,300); 
ylim([5,7]); set(gca,'FontSize',20);

figure(2); clf; 
br = br_stabl(funcs,br,0,1);
[xm,ym] = df_measr(0,br); 
br_splot(br,xm,ym);
ylim([5,7]); set(gca,'FontSize',20);
\end{verbatim}
\end{scriptsize}
\end{tcolorbox}
We now continue the equilibrium branch \texttt{br} and then determine the stability of the continued branch.  
The above commands yield two plots of the equilibrium branch \texttt{br} as shown in Figure \ref{fig:contnsplot}. 
The plot on the left panel of Figure \ref{fig:contnsplot} shows \texttt{br} with additional 300 points. So now \texttt{br} contains a total of 302 points which is achieved by using the function \texttt{br$\_$contn}. 
Meanwhile, the plot on the right panel of Figure \ref{fig:contnsplot} shows the same branch \texttt{br} but with stability information. 
This is obtained using the function \texttt{br$\_$stabl}. 
The stable and unstable parts of the branch are in green and red, respectively, 
while the Hopf bifurcation points are marked with asterisks ($\ast$). 
It is worth noting that in this particular example the stability switches occur at the Hopf bifurcation points. Here, there are five stability switches. \\

The function \texttt{df$\_$measr} gives the default \texttt{xm} and \texttt{ym} for the equilibrium branch \texttt{br}. We can check what these are by typing \texttt{xm} and \texttt{ym} in the Command Window.

\begin{scriptsize}
\begin{verbatim}
>> xm

xm = field: 'parameter'
     subfield: ''
     row: 1
     col: 8
     func: ''

>> ym

ym = field: 'x'
     subfield: ''
     row: 1
     col: 1
     func: ''
\end{verbatim}
\end{scriptsize}
Hence, for the plots in Figure \ref{fig:contnsplot}, the horizontal axis is terms of the time delay parameter $\tau$ while the vertical axis is in terms of the state variable $x(t)$. If you want the vertical axis to be in terms of $y(t)$ instead of $x(t)$, you need to type in \texttt{ym.row=2} in the Command Window to make the desired change. Here, we keep the vertical axis in terms $x(t)$ for future plots.\\

\begin{tcolorbox}
\begin{scriptsize}
\begin{verbatim}
br.method.stability.minimal_real_part = -2;
nunst = GetStability(br);
ind_hopf = find(abs(diff(nunst))==2)
\end{verbatim}
\end{scriptsize}
\end{tcolorbox}
The function \texttt{GetStability} yields \texttt{nunst} which is the number of characteristic roots in the open right-half plane. Hopf bifurcation occurs when a complex conjugate pair of simple characteristic roots crosses the imaginary axis. The above set of commands gives \texttt{ind$\_$hopf} which is the list of points in the equilibrium branch \texttt{br} where Hopf bifucation occurs.
\begin{scriptsize}
\begin{verbatim}
ind_hopf = 30
           110
           142
           251
           252
\end{verbatim}
\end{scriptsize}
The code below gives an animation showing the movement of the characteristic roots on the complex plane as the time delay parameter $\tau$ is varied. 
Pay attention to the value of $\tau$ when a pair of complex conjugate roots crosses the imaginary axis. 
\begin{scriptsize}
\begin{verbatim}
for i = 1:length(br.point); clf;
    figure(33); hold on;
        plot([-0.5 0.5], [0 0], 'b', [0 0], [-5 5], 'b');
        p_splot(br.point(i));
        axis([-0.20 0.20 -1.5 1.5]);
        tau = br.point(i).parameter(8);
        text(0.01,0.2,['\tau = ',num2str(tau,'%2.2f')],'FontSize',32);
        M(i) = getframe(gcf);    
end
\end{verbatim}
\end{scriptsize}
The value of the time delay parameter $\tau$ at the branch points given in \texttt{ind$\_$hopf} are obtained as follows.
\begin{scriptsize}
\begin{verbatim}
>> for i = 1:length(ind_hopf) 
     critical_tau(i) = br.point(ind_hopf(i)).parameter(8);
   end
>> critical_tau'

ans =  1.348598400000001
       5.348598399999998
       6.948598399999993
      12.398598400000049
      12.448598400000050
\end{verbatim}
\end{scriptsize}
Comparing with values in Table \ref{tab:taus}, the first three values above are correct up the first decimal digit while the last two are not. We have to keep in mind that the above values are mere approximations since the maximum stepsize that we assigned for the time delay parameter $\tau$ in \texttt{parbd} is just \texttt{0.05}.\\

\begin{tcolorbox}
\begin{scriptsize}
\begin{verbatim}
[br_hopf1,success] = SetupHopf(funcs,br,ind_hopf(1))
\end{verbatim}
\end{scriptsize}
\end{tcolorbox}
The function \texttt{SetupHopf} allows us to initialize the continuation of Hopf bifurcations. 
Here, we denote by \texttt{br$\_$hopf1} the Hopf branch which as we saw earlier occurs approximately at the point \texttt{ind$\_$hopf(1)} along the equilibrium branch \texttt{br}. The above command yields

\begin{scriptsize}
\begin{verbatim}
br_hopf1 = method: [1x1 struct]
           parameter: [1x1 struct]
           point: [1x1 struct]

success = 1
\end{verbatim}
\end{scriptsize}
which means that our attempt to set-up the Hopf branch \texttt{br$\_$hopf1} is successful. 
At the first branch point of \texttt{br$\_$hopf1}, the value of the time delay parameter $\tau$ can be obtained by typing in the following in the Command Window.
\begin{scriptsize}
\begin{verbatim}
>> br_hopf1.point(1).parameter(8)

ans = 1.379413927096384
\end{verbatim}
\end{scriptsize}
The above value is a correction to the initial guess \texttt{1.348598400000001} obtained from \texttt{br.point(ind$\_$hopf(1)).parameter(8)}. 
Moreover, this corrected value matches the value of $\tau_0^+$ given in Table \ref{tab:taus}.\\ 

\begin{tcolorbox}
\begin{scriptsize}
\begin{verbatim}
[br_hopf2,success] = SetupHopf(funcs,br,ind_hopf(2))
[br_hopf3,success] = SetupHopf(funcs,br,ind_hopf(3))
[br_hopf4,success] = SetupHopf(funcs,br,ind_hopf(4))
[br_hopf5,success] = SetupHopf(funcs,br,ind_hopf(5))
\end{verbatim}
\end{scriptsize}
\end{tcolorbox}
Similarly, the correct value of $\tau$ for the next four Hopf bifurcation points are obtained by setting up Hopf branches. Here, we denote the next four Hopf branches as \texttt{br$\_$hopf2}, \texttt{br$\_$hopf3}, \texttt{br$\_$hopf4}, and \texttt{br$\_$hopf5}. 
The correct values of $\tau$ for the next four Hopf bifurcation points are obtained as follows.
\begin{scriptsize}
\begin{verbatim}
>> [br_hopf2.point(1).parameter(8);...
    br_hopf3.point(1).parameter(8);...
    br_hopf4.point(1).parameter(8);...
    br_hopf5.point(1).parameter(8)]

ans =  5.393140023781609
       6.981037144585396
      12.538835708751554
      12.538835707843552
\end{verbatim}
\end{scriptsize}
\noindent These values match the values in Table \ref{tab:taus} of $\tau_1^-$, $\tau_1^+$, and $\tau_2^-$, but not that of $\tau_2^+$. We remedy this by adding \texttt{'excludefreqs',br$\_$hopf4.point(1).omega} in the previous code which removes undesired eigenvalues from consideration.

\begin{tcolorbox}
\begin{scriptsize}
\begin{verbatim}
[br_hopf5,success] = SetupHopf(funcs,br,ind_hopf(5),...
    'excludefreqs',br_hopf4.point(1).omega)
\end{verbatim}
\end{scriptsize}
\end{tcolorbox}
\noindent The correct value for $\tau_2^+$ is now obtained by typing in the following in the Command Window.
\begin{scriptsize}
\begin{verbatim}
>> br_hopf5.point(1).parameter(8)

ans = 12.582660362074412
\end{verbatim}
\end{scriptsize}
  
Branches of periodic solutions can be obtained from the identified Hopf bifurcations. 
Among these branches of periodic solutions, the branch emanating from the third Hopf bifurcation shows some interesting dynamics. We focus on this branch for the rest of this section.
\begin{tcolorbox}
\begin{scriptsize}
\begin{verbatim}
[br_psol3,success] = SetupPsol(funcs,br,ind_hopf(3))
\end{verbatim}
\end{scriptsize}
\end{tcolorbox}
The function \texttt{SetupPsol} allows us to set-up a branch of periodic solutions. Here, we denote by \texttt{br$\_$psol3} the branch of periodic solutions that emanates from the third Hopf bifurcation which as we saw earlier occurs approximately at the point \texttt{ind$\_$hopf(3)}  along the equilibrium branch \texttt{br}. The above command yields
\begin{scriptsize}
\begin{verbatim}
br_psol3 = method: [1x1 struct]
           parameter: [1x1 struct]
           point: [1x2 struct]

success = 1
\end{verbatim}
\end{scriptsize}
which means that our attempt to set-up the branch of periodic solutions \texttt{br$\_$psol3} is successful. 

\begin{tcolorbox}
\begin{scriptsize}
\begin{verbatim}
br_psol3.method.continuation.plot = 0; 
[br_psol3,s,f,r] = br_contn(funcs,br_psol3,60);
br_psol3 = br_stabl(funcs,br_psol3,0,1);

xm_psol=xm;
ym_psol.field='profile';
ym_psol.subfield='';
ym_psol.row=1;
ym_psol.col='all';
ym_psol.func='max';

figure(3); clf; hold on;
br_splot(br,xm,ym);
br_splot(br_psol3,xm_psol,ym_psol);
axis([6.5 9.5 0 45]); set(gca,'FontSize',20);
\end{verbatim}
\end{scriptsize}
\end{tcolorbox}

\begin{figure}
\centering
\includegraphics[width=0.5\textwidth]{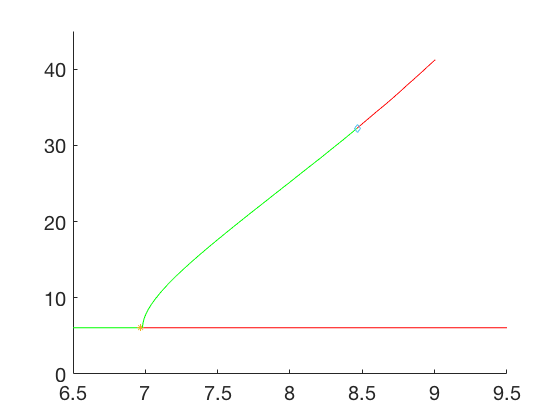}
\caption{Plot of the branch of periodic solutions \texttt{br$\_$psol3} emanating from the third Hopf bifurcation point at $\tau= \tau_2^- = 6.981037144585396$.}
\label{fig:psol3}
\end{figure}

The above commands yield a plot of the periodic-solution branch \texttt{br$\_$psol3} together with the equilibrium branch \texttt{br} as shown in Figure \ref{fig:psol3}. 
Plotting the branch \texttt{br$\_$psol3} requires defining \texttt{xm$\_$psol} and \texttt{ym$\_$psol}. Here, \texttt{xm$\_$psol} is the same as \texttt{xm}, i.e. the time delay parameter $\tau$. 
For \texttt{ym$\_$psol}, we use the maximum function taking just the maximum value of the periodic $x(t)$.

\begin{tcolorbox}
\begin{scriptsize}
\begin{verbatim}
nunst_psol3 = GetStability(br_psol3,'exclude_trivial',true);
ind_pd = find(abs(diff(nunst_psol3))==1,1,'first')
[per2,success] = DoublePsol(funcs,br_psol3,ind_pd)
\end{verbatim}
\end{scriptsize}
\end{tcolorbox}  
Observe that \texttt{br$\_$psol3} is initially stable and then it becomes unstable at a point marked with ($\diamond$). 
This stability switch actually occurs at a period-doubling bifurcation. 
We can get the value of $\tau$ where this period-doubling bifurcation occurs by setting up the branch of period-2 solutions \texttt{per2} using the function \texttt{DoublePsol}. 
The value of $\tau$ where the period-doubling bifurcation occurs is obtained by typing in the following commands in the Command Window.

\begin{scriptsize}
\begin{verbatim}
>> per2.point(1).parameter(8)

ans = 8.464201107682122
\end{verbatim}
\end{scriptsize}

\begin{figure}
\centering
\includegraphics[width=0.5\textwidth]{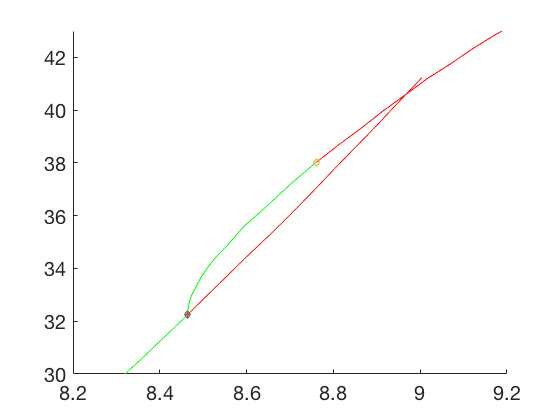}
\caption{Plot of the branch of period-2 solutions \texttt{per2} emanating from the period-doubling bifurcation point ($\diamond$) along the branch of periodic solutions \texttt{br$\_$psol3} where $\tau= 8.464201107682122$.}
\label{fig:per2}
\end{figure}

\begin{tcolorbox}
\begin{scriptsize}
\begin{verbatim}
figure(4); clf; hold on;
per2.method.continuation.plot = 0; 
[per2,s,f,r] = br_contn(funcs,per2,45);
per2 = br_stabl(funcs,per2,0,1);
br_splot(br_psol3,xm_psol,ym_psol);
br_splot(per2,xm_psol,ym_psol);
axis([8.2 9.2 30 43]); set(gca,'FontSize',20);
\end{verbatim}
\end{scriptsize}
\end{tcolorbox}
We continue the branch of period-2 solutions \texttt{per2} emanating from the period-doubling bifurcation  ($\diamond$) along the branch of periodic solutions \texttt{br$\_$psol3}. Figure \ref{fig:per2} shows the plot this continued branch together with the periodic solutions branch \texttt{br$\_$psol3}.

\begin{figure}
\centering
\includegraphics[width=0.3\textwidth]{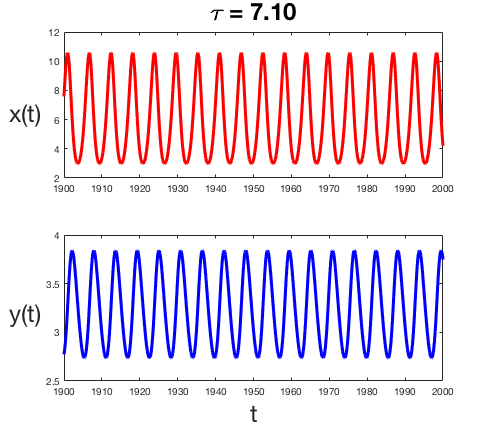}
\includegraphics[width=0.3\textwidth]{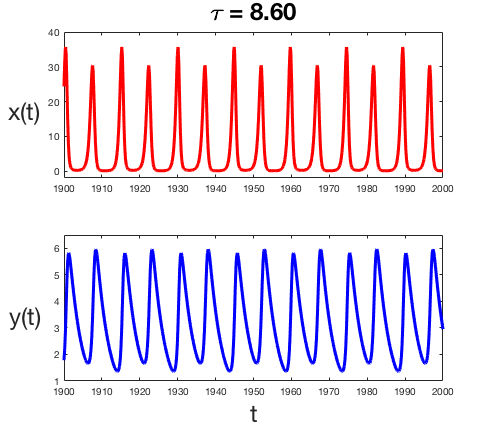}
\includegraphics[width=0.3\textwidth]{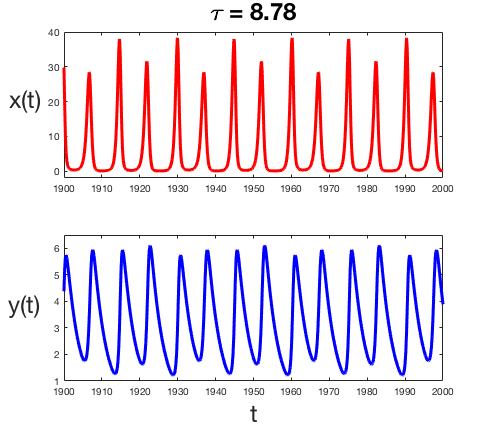}
\caption{Time series plots of $x(t)$ and $y(t)$ for different values of the time delay parameter $\tau$ showing period-1 solutions ($\tau=7.10$), period-2 solutions ($\tau=8.60$), and period-4 solutions ($\tau=8.78$).}
\label{fig:timeseries}
\end{figure}

The period-2 solutions branch \texttt{per2} also undergoes stability switch at a point marked with ($\diamond$). As before, this is also a period-doubling bifurcation. We leave it to the reader to verify that this second period-doubling bifurcation occurs at the value $\tau = 8.757752002502176$. 
That is, beyond this value, we can expect a period-4 solution. Figure \ref{fig:timeseries} shows the time series plots of $x(t)$ and $y(t)$ for different values of the time delay parameter $\tau$ showing period-1 solutions ($\tau=7.10$), period-2 solutions ($\tau=8.60$), and period-4 solutions ($\tau=8.78$).

We end this section with an exercise for the readers to check if the period-doubling bifurcations will go on and will lead to a cascade of period-doubling bifurcations and eventually to chaos.

\section{Conclusions} 
In this paper we presented how to use DDE-Biftool to obtain and analyze branches of solutions to a system of delay differential equations. 
We did this by revisiting the work of Toaha and Hassan \cite{ToahaHassan2008} on a harvested predatory-prey model with time delay. We first discussed this model theoretically using a more simplified approach. Then, we implement the model in DDE-Biftool to study it numerically. The values of the time delay where Hopf bifurcations occur were numerically obtained and matches the theoretical results. In addition, branches of periodic solutions were also obtained using numerical continuation. At one of the branches of periodic solutions, some interesting dynamics occured. The occurrence of period-doubling bifurcations, which could lead to chaotic behavior of the system, was observed. 
Here, we emphasize the importantance of both theory and numerics in studying models. We hope that this paper served its purpose of introducing researchers to time delay systems and its implementation in DDE-Biftool which reveals more dynamical behavior of the model being  considered.

\paragraph{Acknowledgements.}
The author acknowledges the support of University of the Philippines Baguio, CIMPA, IMU-CDC, SEAMS, and Universiti Sains Malaysia for his participation to SEAMS School 2018 on Dynamical Systems and Bifurcation Analysis. 
The author also would like to thank the referees for their valuable reviews that improved the quality of this paper. 

%
%

\end{document}